\documentclass[12pt]{amsart}
\usepackage{geometry}          
\usepackage{graphicx}
\usepackage{amssymb}
\usepackage{epstopdf}
\usepackage{ stmaryrd }
\def\lab#1{\label{#1}}\def\comment#1{}

\newtheorem{theo}[subsection]{Theorem}

\newtheorem{lemm}[subsection]{Lemma}
\newtheorem{defi}[subsection]{Definition}
\newtheorem{prop}[subsection]{Proposition}
\newtheorem{condi}[subsection]{Condition}

\newtheorem{rema}[subsection]{Remark}

\def\equat{\refstepcounter{subsection}$$~}

\def\endequat{\leqno{(\boldsymbol{\arabic{section}.\arabic{subsection}})}~$$}

\def\F{{\Bbb{F}}}
\def\C{{\Bbb{C}}}
\def\N{{\Bbb N}}

\def\rrr{{\bf r}}\def\sss{{\bf s}}\def\ttt{{\bf t}}

\def\LLL{{\mathcal{L}}}
\def\MMM{{\mathcal{M}}}

\def\PPP{{\mathcal{P}}}
\def\QQQ{{\mathcal{Q}}}
\def\RRR{{\mathcal{R}}}
\def\SSS{{\mathcal{S}}}

\def\VVV{{\mathcal{V}}}
\def\WWW{{\mathcal{W}}}

\def\aa{{\alpha}}
\def\bb{{\beta}}
\def\ga{{\gamma}}

\def\DD{{\Delta}}

\def\si{{\sigma}}
\def\cc{{\chi}}

\def\la{{\lambda}}

\def\Om{{\Omega}}
\def\oo{{\omega}}
\def\zz{{\zeta}}

\def\Sy{{\frak S}}

\def\TTTT{{\frak T}}

\def\ti#1{\widetilde{#1}}
\def\ha#1{\widehat{#1}}

\usepackage{amsmath}
\usepackage{accents}

\def\Id{{\rm Id}}

\def\Aut{{\rm Aut}}
\def\GL{{\rm GL}}\def\PSL{{\rm PSL}}\def\PSU{{\rm PSU}}\def\Sp{{\rm Sp}}

\def\Res{{\rm Res}}

\def\Ce#1{{{\rm C}_{#1} }}
\def\Ze#1{{\rm Z}(#1)}
\def\No#1{{{\rm N}_{#1}}}

\def\ad{{\rm ad}}

\def\Irr{{\rm Irr}}\def\IBr{{\rm IBr}}
\def\AWe{{\rm Alp}}

\def\conj{{\bf conj}}

\def\inn{\subseteq}
\def\nni{\supseteq}
\def\tr{\unlhd}
\def\smd{\rtimes}

\def\mm{{-1}}

\def\te{\otimes}

\def\pp{{p'}}

\def\Ct{{\C^\times}}\def\Ft{{\F^\times}}

\def\qed{\hfill\vrule height 1.6ex width .7ex depth -.1ex }

\def\mr#1{\smash{\mathop{\relbar\joinrel\longrightarrow}\limits^{#1}}}

\def\Pr#1{Proposition~\ref{#1}}\def\Th#1{Theorem~\ref{#1}}\def\Le#1{Lemma~\ref{#1}}\def\De#1{Definition~\ref{#1}}

\title[On the reduction of Alperin's weight conjecture]{Two remarks on the reduction of Alperin's weight conjecture.}
\author{ Marc Cabanes}

\begin{document}

\maketitle


\begin{abstract}
The so-called inductive McKay condition on finite simple groups, due to Isaacs-Malle-Navarro, has been recently reformulated by Sp\"ath. We show that this reformulation applies to the reduction theorem for Alperin's weight conjecture, due to Navarro-Tiep. This also simplifies the checking of the inductive condition for Alperin's weight conjecture in the case of simple groups of Lie type with regard to the defining prime.
\end{abstract}

\let\thefootnote\relax

\section*{ Introduction}

Thanks to Isaacs-Malle-Navarro's reduction theorem [IMN], McKay's conjecture on complex characters of finite groups reduces to a so-called inductive condition that must be fulfilled by any finite quasi-simple group. This has opened the way to a proof of McKay's conjecture using the classification of those finite groups. Isaacs-Malle-Navarro's inductive condition requires of course that each quasi-simple group satisfies McKay's conjecture, but it also requires the fulfilment of another condition, of a cohomological nature (see [IMN]~\S 10 (8)). The whole inductive condition has been reformulated by Sp\"ath in [S1], [S2]. This brings a major simplification to both the proof of the reduction theorem and to the actual verification of the cohomological condition, allowing it in particular for quasi-simple groups of Lie type with regard to the defining prime [S2].

Let us recall now Alperin's weight conjecture [Al]. When $p$ is a prime and $G$ is a finite group, let $\F$ be an algebraically closed field of characteristic $p$ and $\F G$ the corresponding group algebra. Denote by $\IBr (G)$ the set of simple $\F G$-modules up to isomorphism, and by $\AWe (G)$ the set of $G$-conjugacy classes of pairs $(P,\psi )$ where $P$ is a $p$-subgroup of $G$ and $\psi\in \IBr (\No G(P)/P)$ is a simple \it projective \rm $\F (\No G(P)/P)$-module. Alperin conjectured that $$|\IBr (G)|=|\AWe (G)|.$$

Navarro-Tiep [NvTi] have shown that this conjecture is implied by a related ``inductive" condition to be fulfilled by any finite quasi-simple group. We show here that the main ideas of [S1], [S2] apply also to Brauer characters and the crucial step of Navarro-Tiep's reduction in [NvTi]~3.2. This reduction step is the one common to Alperin's weight conjecture and McKay's conjecture (see [IMN]~13.1).

Our main result (\Th{MainT}) shows that Sp\"ath's version of the cohomological condition used in [IMN] and [NvTi] implies an isomorphism of certain character triples. In the proof of the reduction theorems, it is through this statement that the condition on simple subquotients of an arbitrary group $G$ is used.

We also introduce a condition on modular representations ensuring this cohomological condition (see \Pr{MQ}). This simplifies the checking of the inductive conduction for Alperin's weight conjecture for quasi-simple groups of Lie type with regard to the defining prime, which is known from [NvTi] Theorem~C.

\section{Notations and background results.}

\noindent \bf 1.1. Groups. \rm 
Assume the group $G$ acts on the set $X$. If $x\in X$, one denotes $G_x:=\{ g\in G\mid g.x=x\}$ the stabilizer of $x$. If $G$ acts on several sets $X$, $X'$,$\dots$ and $x\in X$, $x'\in X'$, $\dots$, then $G_{x,x',\dots }$ denotes the intersection of stabilizers of $x$, $x'$, $\dots$.

For $G$ a group, its automorphism group is denoted by $\Aut (G)$. If $g\in G$, $\ad_g\in\Aut (G)$ denotes the inner automorphism $x\mapsto \ad_g(x)={}^gx=gxg^\mm$. 

When $H$ is a subgroup of $G$, one denotes by $\No G(H)\nni \Ce G(H)$, its normalizer, resp. centralizer in $G$. A (normalized) $H$-section $$\sss\colon G/H\to G$$ is any abstract map such that $\sss (1)=1$ and $\sss (x)H=x$ for any $x\in G/H$. If $H,K$ are subgroups of $G$, one denotes by $[H,K]$ the subgroup generated by commutators $[h,k]=hkh^\mm k^\mm$ for $h\in H$, $k\in K$. One calls $G$ perfect whenever $G=[G,G]$. A quasi-simple group is a perfect group $G$ such that $G/\Ze G$ is simple.

If $p$ is a prime, one denotes by $G_p$, $G_\pp$, and ${\rm O}_p(G)$ $\dots$ the subset of $p$-elements, resp. $p$-regular elements, and maximal normal $p$-subgroup of $G$.

\medskip

\noindent \bf 1.2. Modular representations. \rm (See [Nv]~\S 8.) \rm
Let $\F$ be an algebraically closed field. If $G$ is a finite group, we denote by $\IBr (G) =\Irr (\F G)$ the set isomorphism classes of irreducible $F$-linear representations of $G$ or equivalently of simple $\F G$-modules. If $H$ is a subgroup and $\zeta\in \IBr (H)$, one denotes by
 $\IBr (G|\zeta )$ the subset of $\IBr ( G)$ corresponding with simple $\F G$-modules whose restriction to $\F H$ contains a submodule in the isomorphism class of $\zz$. We won't really use the interpretation in terms of characters, except for the notation $\zz (1)\in\N $ to denote the dimension of the representation involved, along with the case of $1$-dimensional representations which we identify with group morphisms $G\to\Ft :=\GL_1(\F)$. An irreducible linear representation $\LLL\colon G\to\GL_n(\F)$ (or a  simple $\F G$-module) is said to afford $\cc\in\IBr (G)$ if $\cc$ is its isomorphism class.
 
\medskip

\noindent \bf 1.3. Cocycles and Schur multipliers. \rm (See [As]~\S 33, [NgTs]~\S 3.5.)
If $A$ is a commutative (multiplicative) group and $G$ a finite group, let us recall the (multiplicative) groups $A ^G$ of arbitrary maps $G\to A$, and ${\rm Z}^2 (G,A )$ of cocycles $G\times G\to A$. Recall $$\partial \colon {A }^G\to {\rm Z}^2 (G, A)\ ,\ (\partial f)(g_1,g_2)=f(g_1)f(g_2)f(g_1g_2)^\mm$$ for any $f\colon G\to A$ and $g_1,g_2\in G$. Recall ${\rm H}^2 (G, A )={\rm Z}^2(G, A)/\partial ({A}^G)$.

Recall the Schur multiplier ${\rm H}^2 (G,\Ct)$, which is a finite group. A covering of a \it perfect \rm finite group $G$ is a surjective group morphism $H\to G$ with perfect $H$ and central kernel. A universal covering $\ha G\to G$ is one such that its kernel is the Schur multiplier of $G$. It is unique up to isomorphism of coverings.

\medskip

\noindent \bf 1.4. Projective representations. \rm (See [CuRe]~\S 11, [NgTs]~\S 3.5, [Nv]~\S 8.)
A (normalized) projective representation of a group $G$ is a map $\PPP\colon G\to \GL_n(\F )$ such that $\PPP (1)=\Id_n$ and there is a (unique) $\aa_\PPP\in {\rm Z}^2(G,\Ft )$ such that $\PPP (g)\PPP (g')=\aa_\PPP (g,g').\PPP (gg')$ for any $g,g'\in G$. The integer $n$ is called the dimension of $\PPP$. Moreover $\PPP$ is called irreducible if and only if no proper subspace of $\F^n$ is stable by $\PPP (G)$. Note that $\Aut (G)$ acts on projective representations of $G$ by $^\tau\PPP:=\PPP\circ\tau^\mm$.

The linear representations are the projective representations whose cocycle is trivial (identically $1$). 
Changing a projective representation $\PPP$ into a proportional representation $\bb \PPP $ with $\bb \in\Ft^G$ will change the associated cocycle by $\partial\bb$, namely $\aa_{\bb\PPP}=\partial\bb .\aa_\PPP$.

 Recall that two projective representations $\mathcal{X}$, $\mathcal{X}'\colon X\to\mathrm{GL}_n(\F )$ of a finite group $X$ are said to be \it equivalent \rm if and only if there is $J\in \mathrm{GL}_n(\F )$ such that $$J\mathcal{X}'(x)J^{-1}\in\Ft \mathcal{X}(x)$$ for any $x\in X$. 
We will need the following elementary facts.

 \begin{prop}\label{ExtPro}  {\sl Let $X\tr G$ be finite groups and let $\VVV$ be an irreducible projective representation of $X$ stabilized by $G$ (that is $\VVV\circ\ad_g$ is equivalent to $\VVV $ for any $g\in G$). Assume $X$ is perfect. Then 
 
 (i) $\VVV$ extends into a projective representation $\ti\VVV$ of $G$, unique up to equivalence.
 
 (ii) If moreover the cocycle of $\VVV$ is of finite order (in ${\rm Z}^2(X,\Ft)$), one may choose $\ti\VVV$ so that its cocycle has finite order (in ${\rm Z}^2(G,\Ft)$).
 }
\end{prop} 

\noindent{\it Proof of~\Pr{ExtPro}.}
(i) Arguing as in the proof of [NgTs]~3.5.7.(i) or [Nv] 8.14, one finds a map $\widetilde{\mathcal{V}}\colon G\to \mathrm{GL}_{d}(\F)$ extending $\mathcal{V}$ and such that, fixing $a,a'\in G$ and defining $T=\widetilde{\mathcal{V}}(a)\widetilde{\mathcal{V}}(a')\widetilde{\mathcal{V}}(aa')^{-1}$, one has $$\mathcal{V}(x)T\mathcal{V} (x)^{-1}\in\Ft  .T\ \ {\rm for\ all}\ x\in X.$$ The invertible scalar above is a multiplicative map $X\to\Ft$, while $X$ is perfect. So the above equation actually is a commutation $\mathcal{V}(x)T=T\mathcal{V}(x)$ for all $x\in X$, and Schur's lemma tells us that $T$ is scalar. So $\widetilde{\mathcal{V}}$ is a projective representation.

If $\WWW$ is another projective representation of $G$ extending $\VVV$, one has the same equation as above with $T:=\ti\VVV (g)\WWW (g)^\mm$. So the same argument as above implies that $\ti\VVV$ and $\WWW$ are proportional, hence equivalent.

(ii) The group ${\rm H}^2(G,\Ft)$ is finite, so we may choose an integer $N\geq 1$ such that $(\aa_\VVV)^N =1$ in ${\rm Z}^2(X,\Ft)$ and $h^N=1$ for any $h\in {\rm H}^2(G,\Ft)$. If the characteristic of $\F$ is non zero, one can take $N$ coprime to it, so that any element of $\Ft$ has an $N$-th root. 

We have $(\aa_{\ti\VVV})^N\in \partial ((\Ft)^G)$, so $(\aa_{\ti\VVV})^N=\partial\ga$ where $\ga\in (\Ft)^G$. Note that from our hypothesis and the choice of $N$, $\partial\ga$ is trivial on $X$. This means that the restriction of $\ga$ to $X$ is a multiplicative map $X\to\Ft$. But this can only be trivial, so $\ga (X)=\{ 1\}$. Let's now take an $N$-th root $\bb (g)$ of each $\ga (g)^\mm$ ($g\in G$) with $\bb (X)=\{ 1\}$. Then $\bb\ti\VVV$ coincides with $\ti\VVV$, hence $\VVV$ on $X$, and its cocycle is $\partial\bb.\aa_{\ti\VVV}$ whose $N$-th power is $\partial (\bb^N)\partial (\ga)=\partial (\bb^N\ga)=1$.

\qed

\section{Embeddings and cohomologically related representations.}

We fix an algebraically closed field $\F$. All simple modules for group algebras are over $\F$. We use the notation $\IBr (G)$ (which may of course mean $\Irr (G)$ when the characteristic of $\F$ is 0).

One recalls Sp\"ath's condition on pairs of linear irreducible representations (see [S1]~2.8). In the following, $ X\nni Y\nni Z = \Ze X$ are finite groups and $\cc\in \IBr (X)$, $\psi\in \IBr (Y)$.

\begin{condi}\label{Coho} \rm {\rm
There are (normalized) projective representations $\PPP$ of $\Aut (X)_\cc$, $\QQQ$ of $\Aut (X)_{Y,\psi}$ such that (identifying $X/Z$ with the subgroup of inner automorphisms in $\Aut (X)$)

(i) $\Res^{\Aut (X)_\cc}_{X/Z}\PPP =\LLL\circ \sss$ where $\sss\colon X/Z\to X$ is a (normalized) section of the quotient map $X\to X/Z$, and $\LLL$ is a linear representation of $X$ affording $\cc$,

(i') $\Res^{\Aut (X)_{Y,\psi}}_{Y/Z}\QQQ=\MMM\circ \ttt $ where $\ttt\colon Y/Z\to Y$ is the restriction of $\sss$, and $\MMM$ is a linear representation of $Y$ affording $\psi$,

(ii) $\Res^{\Aut (X)_\cc}_{{\Aut (X)_{Y,\psi,\cc}}}\PPP$ and $\Res^{\Aut (X)_{Y,\psi}}_{{\Aut (X)_{Y,\psi,\cc}}} {\QQQ}$ have same cocycle $\ga\in {\rm Z}^2({\Aut (X)_{Y,\psi,\cc}},\Ft)$.

}\end{condi}

\begin{rema}\label{CohoRk} {\rm  (i) Let $(\cc ,\psi )\in \IBr (X)\times \IBr (Y)$ and $\si\in\Aut (X)$. Then $(\cc ,\psi )$ satisfies Condition~\ref{Coho} if and only if $(\cc^\si ,\psi{}^\si )\in \IBr (X)\times \IBr (Y{}^\si )$ satisfies Condition~\ref{Coho}.

(ii)
Take $\LLL$ an irreducible linear representation of $X$. If $\sss\colon X/Z\to X$ is a section for $Z=\Ze X$, one may define $\partial\sss\in {\rm Z}^2(X/Z,Z)$ by $$\partial\sss (a,a')=\sss (a)\sss (a')\sss (aa')^\mm$$ for $a,a'\in X/Z$. Then denoting by $\la\colon Z\to \Ft$ the linear character such that $\LLL (z)=\la (z)\Id_{\cc (1)}$, one has clearly $\aa_{\LLL\circ \sss}=\la\circ\partial\sss$. This is of finite order. 

So, when $X$ is perfect, \Pr{ExtPro} applied to $X/Z\tr\Aut (X)_\cc$ tells us that $\PPP$ as in Condition~\ref{Coho} always exists and has a finite cocycle.

}\end{rema}

Here are some less elementary properties.

\begin{prop}\label{CohQu}  {\sl  Assume $ X\nni Y\nni Z = \Ze X$ are finite groups, and that
$(\cc,\psi)\in \IBr (X)\times \IBr (Y)$ satisfies Condition~\ref{Coho}. Assume $X$ is perfect with trivial Schur multiplier (or just ${\rm H}^2(X,\Ft)=\{ 1\}$). Assume $Z'\inn Z$ is in both kernels of $\cc$ and $\psi$, and denote by $(\cc ',\psi ')\in \IBr (X/Z')\times \IBr (Y/Z')$ the corresponding characters obtained by deflation. Then $(\cc ',\psi ' )$ satisfies Condition~\ref{Coho} for $X/Z'\nni Y/Z'\nni Z/Z' =\Ze {X/Z'}$.}
\end{prop}

\noindent{\it Proof of~Proposition~\ref{CohQu}.} If $X$ is perfect with trivial multiplier, $\Aut (X)_{Z'}$ identifies with $\Aut (X/Z')$ by the universal property of coverings. Then $\Aut (X/Z')$ injects in $\Aut (X)$ by a map which is the identity on inner automorphisms. So it is easy to deduce Condition~\ref{Coho} for  $(\cc ',\psi ' )$ from the same property of $(\cc ,\psi)$ simply by restricting the projective representations $\PPP$ and $\QQQ$ from $\Aut (X)_\cc$ and $\Aut (X)_{Y,\psi}$ to $\Aut (X/Z')_\cc$ and $\Aut (X/Z')_{Y/Z',\psi}$. This gives \Pr{CohQu}.

\qed

\begin{prop}\label{CohPr}  {\sl  Let $n\geq 1$ and $ X\nni Y_i\nni Z = \Ze X$ some finite groups for $i=1,\dots ,n$, with $X$ quasisimple non-abelian.

Assume $(\cc_i,\psi_i)\in \IBr (X)\times \IBr (Y_i)$ for $i=1,\dots ,n$ satisfy Condition~\ref{Coho}. 

Assume that, 

\noindent{\rm (*)} for any $i,j\in \{ 1,\dots ,n\}$, $\tau\in\Aut (X)$, we have $\cc_i={}^\tau\cc_j$ if and only if $(Y_i,\psi_i)$ and $^\tau (Y_j,\psi_j)$ are $X$-conjugate.

Then $(\cc_1\te\dots \te\cc_n, \psi_1\te\dots\te\psi_n)$ satisfies Condition~\ref{Coho} for $X^n\nni Y_1\times\dots\times Y_n\nni Z^n =\Ze {X^n}$.}
\end{prop}

\noindent{\it Proof of~Proposition \ref{CohPr} \rm (see [S2]~3.4).} Using Remark~\ref{CohoRk}.(i) with elements of $\Aut(X)^n\inn\Aut (X^n)$, we may assume that for any $i,j$, $$\cc_i\in\Aut (X).\cc_j \text{ if and only if }\ \cc_i=\cc_j.$$  Thanks to assumption (*), applying possibly inner automorphisms of $X$ (which leave the $\cc_i$'s unchanged), we may assume that for any $i,j$

(**) $\cc_i\in\Aut (X).\cc_j \Leftrightarrow \cc_i=\cc_j\Leftrightarrow (Y_i,\psi_i)=(Y_j,\psi_j)\Leftrightarrow (Y_i,\psi_i)\in\Aut (X).(Y_j,\psi_j)$. 

On the other hand, since $X$ is quasi-simple non-abelian, each automorphism of $X^n$ has to permute the given summands which are the minimal normal perfect subgroups. So, denoting by $\Sy_n$ the symmetric group on $n$ letters, one has $\Aut (X^n)=\Aut (X){\wr}\Sy_n$. 

Denote $Y=Y_1\times\dots\times Y_n$, $\cc=\cc_1\times\dots\times \cc_n$, $\psi =\psi_1\times\dots\times \psi_n$. With (**) above, it is easy to see that both $\Aut (X^n)_{\cc}$ and $\Aut (X^n)_{Y,\psi}$ split as a direct product along the partition of $\{1,\dots ,n\}$ that (**)  defines. So we may assume that $\cc_i=\cc_1$, $(Y_i,\psi_i)=(Y_1,\psi_1)$ for all $i$'s.

 The fact that $(\cc_1 ,\psi_1)$ satisfies Condition~\ref{Coho} ensures the existence of a (normalized) section $\rrr\colon X/Z\to X$ and projective representations $Ê\PPP_1\colon A\to\GL_\F(V)$ of $A:=\Aut (X)_{\cc_1}$ extending $\LLL_1\circ\rrr$ for $\LLL_1$ a linear representation of $X$ affording $\cc_1$, and $\QQQ_1\colon B\to\GL_\F(W)$ of $B:=\Aut (X)_{Y_1,\psi_1}$ extending $\MMM_1\circ\rrr_{| Y_1/Z}$ for $\MMM_1$ a linear representation of $Y$ affording $\psi_1$ and such that their cocycles coincide on $A\cap B$. Then one may choose $\rrr^n\colon (X/Z)^n\to X^n$ as section and define $\widetilde{\PPP}_1\colon A{\wr}\Sy_n\to\GL_\F (V^{\te n})$ on $\Aut (X^n)_\cc =A{\wr}\Sy_n$ by the usual construction of linear representations of wreath products (see [B]~\S 3.15) suitably generalized to projective representations. Namely, for $a_1,\dots ,a_n\in A$, $\si\in\Sy_n$, $$\widetilde{\PPP}_1\bigl( (a_1,\dots ,a_n).\si\bigr)=(\PPP_1(a_1)\te\dots\te\PPP_1(a_n))\circ \si_V$$ where $\si_V$ is the endomorphism of $V^{\te n}$ sending $v_1\te\cdots \te v_n $ to  $v_{\si^\mm (1)}\te\cdots \te v_{\si^\mm (n)}$. As for cocycles, we have $\aa_{\widetilde{\PPP}_1}\big( (a_1,\dots ,a_n),\si ;  (a'_1,\dots ,a'_n),\si '\big)=$\newline\rightline{$\aa_{\PPP_1} (a_1,a'_{\si^{-1}(1)})\dots \aa_{\PPP_1} (a_n,a'_{\si^{-1}(n)})$. }

The same process yields some $\widetilde{\QQQ}_1$ on $B{\wr}\Sy_n=\Aut (X^n)_{Y,\psi}$ from $\QQQ_1$ with the corresponding properties and formula for its cocycle. It is then clear that $\aa_{\widetilde{\PPP}_1}$ and $\aa_{\widetilde{\QQQ}_1}$ have same restriction to $(A\cap B){\wr}\Sy_n$ since $\aa_{\QQQ_1}$ and $\aa_{\PPP_1}$ have same restriction to $A\cap B$.

This completes the checking of Condition~\ref{Coho}.

\qed

 \begin{theo}\label{MainT}  {\sl Let $Z=\Ze X\inn Y\inn X\tr G$ be finite groups,  and let $\nu \in\IBr (Z)$, $\cc \in\IBr (X | \nu)$, $\psi \in\IBr (Y |\nu)$.
 Denote $X'=[X,X]$, $Y'=Y\cap X'$, $Z'=Z\cap X'$.
 
 Assume the following hypotheses 
 
 (i) $X=X'Z$ (so that $X'$ is perfect with center $Z'$) and $Z\inn \Ze G$.
 
 (ii) $\No X(Y)_\psi =Y$, $\No G(Y)_\cc =\No G(Y)_\psi$ and $G_\cc =X.\No G(Y)_\psi $ (so that $ \No G(Y)_\psi /Y \cong G_\cc /X$ by the natural map).
 
(iii) $\Res^X_{X'}\cc $ and $\Res^Y_{Y'}\psi $ satisfy Condition~\ref{Coho} with respect to $X'\nni Y'\nni Z'$. 

\smallskip

Then the modular character triples $(G_\cc , X,\cc)$ and $(\No G(Y)_\psi ,Y,\psi)$ are isomorphic in the sense of [Nv]~8.25. In particular $$| \IBr(G |\cc )|=|\IBr (\No G(Y) | \psi )|.$$ }
\end{theo}

\noindent{\it Proof of~\Th{MainT}.}  As said before, most arguments are taken from [S2]~\S 3.

Assume first that $G_\cc = G$, and therefore $\No G(Y)_\psi =\No G(Y)$. 

By \De{Coho}, there are normalized projective representations $\PPP '$, $\QQQ '$ of $\Aut (X')_{\Res^X_{X'}\cc }$, resp. 
$\Aut (X')_{Y', \Res^Y_{Y'}\psi }$ such that their restrictions to $X'/Z'$, resp. $Y'/Z'$ are $\LLL '\circ \sss$, $\MMM '\circ\ttt$ where $\LLL '$, $\MMM '$ are linear representations affording $ \Res^X_{X'}\cc $, $\Res^Y_{Y'}\psi $ with $\sss\colon X'/Z'\to X'$ a section, and $\ttt$ its restriction to $Y'/Z'$, with moreover $\aa_{\PPP '}$ and $\aa_{\QQQ '}$ coinciding on $\Aut (X')_{Y', \Res^Y_{Y'}\psi ,\Res^X_{X'}\cc }$.

We have $X'/Z'\inn G/\Ce G(X')=G/\Ce G(X)\inn \Aut (X')_{\Res^X_{X'}\cc }$ so we may restrict the above projective representations to the corresponding intermediate groups $G/\Ce G(X)$ and $\No G(Y)/\Ce G(X)$. We may notice also that $\LLL '\circ \sss =\LLL\circ \sss$ and $\MMM '\circ\ttt = \MMM \circ\ttt$ where $\sss$, resp. $\ttt$ is considered having values in $X$, resp.  $Y$ and $\LLL$, resp. $\MMM$, is a linear representation affording $\cc$, resp. $\psi$. 

To sum up, we get projective representations $\PPP$, $\QQQ$ of $G/\Ce G(X)$ , $\No G(Y)/\Ce G(X)$ such that their restrictions to $X/Z$, $Y/Z$ are $\LLL\circ \sss$, $\MMM\circ \ttt$ and $\aa_\QQQ $ is the restriction of $\aa_\PPP$ to $\No G(Y)/\Ce G(X)\times \No G(Y)/\Ce G(X)$.

Let us extend $\sss\colon X/Z=X'/Z'\to X'$ into a $\Ce G(X)$-section $\rrr\colon G/\Ce G(X)\to G$.
Let $\mu\colon \Ce G(X)\to\Ft $ be any map extending $\nu\colon Z\to \Ft$ and such that $\mu (cz)=\mu (c)\nu (z)$ for any $z\in Z$, $c\in\Ce G(X)$.

Define $\PPP_\cc$ a map on $G$ by $$\PPP_\cc (c.\rrr (a))=\mu (c)\PPP (a)$$ for $a\in  G/\Ce G(X)$, $c\in\Ce G(X)$. 
Analogously, one defines $\QQQ_{\psi}$ on $\No G(Y)$ from $\QQQ$ by $$\QQQ_{\psi}(c.\rrr (a))=\mu (c)\QQQ(a)$$ for $a\in \No G(Y)/\Ce G(X)$, $c\in\Ce G(X)$.

\begin{lemm}\label{Ptheta} {\sl  (i) $\PPP_\cc$ is a normalized projective representation of $G$ whose restriction to $X$ is a linear representation affording $\cc$.

(i') $\QQQ_{\psi}$ is a normalized projective representation whose restriction to $Y$ is a linear representation affording $\psi$. 

(ii) $\PPP_\cc (txt^{-1})={\PPP_\cc (t)}\PPP_\cc (x){\PPP_\cc (t)}^\mm$ for any $x\in X$, $t\in \No G(Y)$ ; or equivalently~:

(ii') $ \aa_ {\PPP _{\cc}}(t,x) =\aa_ {\PPP _{\cc}}(txt^\mm ,t) $ for any $x\in X$, $t\in \No G(Y)$.

(iii) $\aa_{\QQQ_{\psi}}$ is the restriction to $\No G(Y)\times \No G(Y)$ of $\aa_{\PPP _{\cc}}$}\end{lemm}

Let us say how this will complete the proof of \Th{MainT}. Let $\TTTT$ be a representative system of $\No G(Y)/Y$ with $1\in\TTTT$. By the hypothesis (ii) of the theorem, this is also a representative system for $G/X$. 

Let $\RRR$, resp. $\SSS$, defined on $G$, resp. $\No G(Y)$, by $$\RRR (xt)=\PPP_\cc (x)\PPP_\cc (t)=\aa_{\PPP_\cc}(x,t)\PPP_\cc (xt)$$ {\rm and} $$ \SSS (yt) =\QQQ_{\psi} (y)\QQQ_{\psi} (t)=\aa_{\QQQ_\psi}(x,t)\QQQ_\psi (yt)$$ for $x\in X$, $y\in Y$, $t\in \TTTT$. 

By \Le{Ptheta}.(i) and (i'), $\Res^G_X\RRR$ and $\Res^{\No G(Y)}_Y\SSS$ are linear representations affording $\cc$ and $\psi$. By the above definition, the maps $\RRR$ and $\SSS$ are proportional to $\PPP_\cc$, resp. $\QQQ_\psi$, hence are projective representations with cocycles $\aa_\RRR=\aa_{\PPP_\cc}.\partial \bb_{ \cc}$, resp. $\aa_{\SSS }=
\aa_{\QQQ_{\psi}}.\partial \bb_{ {\psi}}$ where $\bb_{ \cc}\colon G\to \Ft $, resp. $\bb_{ \psi}\colon \No G(Y)\to \Ft$, are defined by $\bb_{ \cc}(xt)=\aa_{\PPP_\cc}(x,t)$, resp.  $\bb_{ {\psi}}(yt)=\aa_{\QQQ_{\psi}}(y,t)$, for $t\in \TTTT$, $x\in X$, $y\in Y$. 

So Lemma~\ref{Ptheta}.(iii) gives that $\aa_{\SSS }$ is the restriction of $\aa_\RRR$ to $\No G(Y)\times \No G(Y)$. 

By Lemma~\ref{Ptheta}.(ii) above, we have \equat\lab{RRR}  \RRR (xt)=\RRR (x)\RRR (t) \ \ {\rm and}\ \ \RRR (tx)=\RRR ((txt^\mm )t)=\RRR (t)\RRR (x)\endequat for any $x\in X$, $t\in \TTTT$. By \Le{Ptheta}.(ii') and (iii) we also have ${\QQQ_{\psi}}(tyt^\mm)={}^{\QQQ_{\psi}(t)}{\QQQ_{\psi}}(y)$ for any $t\in\TTTT$ and $y\in Y$, and consequently as above, \equat\lab{SSS}  \SSS (yt)=\SSS (y)\SSS (t)\ \ {\rm and}\ \ \SSS (ty)=\SSS (t)\SSS (y).\endequat Since $\RRR$ and $\SSS$ are linear on $X$ and $Y$ respectively, it is easily checked that (\ref{RRR}) holds with any $t\in X$ and (\ref{SSS}) holds with any $t\in\No G(Y)$.

This is precisely the requirement that the projective representations $\RRR$ and $\SSS$ have to satisfy (see [NgTs]~3.5.7, [Nv]~8.14) to determine the isomorphism types of the modular character triples $(G,X,\cc )$ and $(\No G(Y),Y,\psi)$ via a standard application of Clifford theorems ([CuRe]~11.20, [NgTs]~3.5.8, [Nv]~8.16). Moreover we have seen that the cocycles $\aa_\RRR$ and $\aa_\SSS$ identify with the same cocycle on $G/X\cong \No G(Y)/Y$. Hence the isomorphism of modular triples.

All this was in the case $G_\cc = G$. In the general case, applying the above to $G_\cc$,
one gets that the triples $(G_\cc ,X,\cc)$ and $(\No  {G_\cc} (Y),Y,\psi )=(\No G(Y)_\psi ,Y,\psi )    $ are isomorphic. Then Clifford correspondence (see [NgTs]~3.3.2) implies $$|\IBr (G | \cc )| = |\IBr (G_\cc | \cc )| = |\IBr (\No G(Y)_\psi  | \psi )| = |\IBr (\No G(Y) | \psi )| .$$

\qed

\medskip \noindent{\it Proof of Lemma~\ref{Ptheta}.} Let us denote $C:=\Ce G(X)$ and $\pi\colon G\to G/C$ the quotient map.

(i) Writing $\PPP_\chi (g)=\mu (g.\rrr(gC)^\mm)\PPP (gC)$ for any $g\in G$, the map $\PPP_\cc$ is clearly proportional to $\PPP\circ \pi$, so it is a projective representation of dimension $\cc (1)$. Let us now look at its restriction to $X$. In order to get our claim, it suffices to show that for any $z\in Z$, $x'\in X'$

\noindent (a)  $\PPP_\cc (z) =\nu (z)\Id_{\cc (1)}$,

\noindent (b) $\PPP_\cc (x')=\LLL (x')$,

\noindent (c) $\PPP_\cc (zx')=\PPP_\cc (z)\PPP_\cc (x')$.

When $x\in X$, $xC\in X'C$, so $\rrr (xC)=\sss (xC)\in X'$ and $x.\sss (xC)^{-1}\in Z$, so that $$\PPP_\cc (x)=\nu (x.\sss (xC)^{-1})\LLL (\sss (xC)).$$

This implies (a) since $\sss (1)=1$. This also implies that $\PPP_\cc$ is normalized. We have (b) since more generally $\LLL (z'x')=\nu (z')\LLL (x')$ for any $x'\in X'$, $z'\in Z\cap X'$ and $\sss (x'C)\in x'Z'$. Then (c) is clear from the above.

The proof of (i') follows the same lines replacing $G$ with $\No G(Y)$, $X'$ with $Y'$ and $X$ with $\No G(Y)$.

(ii) By the definition of $\aa_ {\PPP _{\cc}}$, one has   $${\PPP_\cc (t)}\PPP_\cc (x){\PPP_\cc (t)}^\mm =\omega_t(x).\PPP_\cc (txt^{-1})$$ for $\omega_{t}(x):=\aa_ {\PPP _{\cc}}(t,x)\aa_ {\PPP _{\cc}}(^tx,t)^{-1}\in\F^\times$. Since $\PPP_\cc$ is linear on $X$, one gets that $\omega_{t}$ is a group morphism $X\to\F^\times$. But the equality $X=[X,X].Z$ implies that $\omega_t$ is defined by its restriction to $Z$. On the other hand, (ii) is clear for $x\in Z$ since then $\PPP_\cc (x)$ is the scalar matrix $\nu (x)\Id_{\cc (1)}$ and $txt^{-1}=x$. So $\omega_t = 1$ for any $t\in \No G(Y)$. Hence (ii) and (ii').

(iii) By the definition of ${\PPP _{\cc}}$ from $\PPP$, its cocycle is $\partial\bb_\mu . \aa_ {\PPP \circ \pi}=\partial\bb_\mu . \aa_ {\PPP} \circ \pi$ where $\bb_\mu\colon G\to\F^\times$ is $g\mapsto \mu (g.\rrr(gC)^{-1})$. Restricting to $\No G(Y)$ will give the cocycle of ${\QQQ_{\psi}}$ since the restriction of $\aa_\PPP$ to ${\No G(Y)/C}Ê\times {\No G(Y)/C}$ is $\aa_{\QQQ}$ by hypothesis.
\qed

\section{Application to Alperin's weight conjecture.}

Let's fix a prime $p$, and an algebraically closed field $\F$ of characteristic $p$. 

If $Q\inn G$ is a $p$-subgroup of the finite group $G$, one denotes by $\IBr_Q(\No G(Q))$ the set of (isomorphism classes) of simple $\F\No G(Q)$-modules which are projective modules when seen as $\F (\No G(Q)/Q)$-modules. This makes sense since simple $\F\No G(Q)$-modules must have ${\rm O}_p(\No G(Q))$, hence $Q$, in their kernel. Note that $\IBr_Q(\No G(Q))\not=\emptyset$ then implies that $Q$ is \it $p$-radical\rm , i.e. ${\rm O}_p(\No G(Q))=Q$. Of course, $\IBr_Q(\No G(Q))^\tau =\IBr_{Q^\tau}(\No G(Q^\tau ))$ for any $\tau\in\Aut (G)$.

Let us denote $$\AWe (G):=\bigl(\coprod_Q\IBr_Q(\No G(Q)\bigr)/G-\conj $$ where the sum is over $p$-subgroups $Q\inn G$ and $G-\conj$ means $G$-conjugacy. One then considers the elements of $\AWe (G)$ as $G$-conjugacy classes of pairs $(Q,\psi )$, the so-called ``weights", with $Q$ a $p$-subgroup of $G$ and $\psi\in\IBr_Q(\No G(Q))$. J. Alperin's ``weight conjecture" (see [Al]~\S 1) is as follows : there is a bijection $$\Omega\colon \IBr (G)\mr{\sim} \AWe(G).\leqno({\rm AWC})$$ The group $\Aut (G)$ acts on $\AWe (G)$ and we also have $\AWe (G)=\coprod_{\nu\in\IBr (\Ze G)}\AWe (G | \nu )$ where $\AWe (G | \nu )$ collects the intersections $\IBr_Q(\No G(Q))\cap\IBr (\No G(Q) |\nu)$ for a given $\nu$. A reasonable requirement is that $\Omega (\IBr (G | \nu))=\AWe (G |\nu )$ for every $\nu$. In [Al], Alperin gives a version of his conjecture where the above $\nu$'s are replaced by blocks of $\F G$, see also a related reduction statement in [P]. 

\it Added in proof: \rm Since submission of the present paper, B. Sp\"ath has completed a reduction theorem for this blockwise version [S3].

 \begin{defi}\label{iAw}  {\rm A group $G$ is said to satisfy {\rm (iAw)} if and only if there is a bijection $\Omega$ as above which is $\Aut (G)$-equivariant, preserves the partition along $\IBr (\Ze G)$, and if for every $\cc\in\IBr (G)$, $(Q,\psi )\in \Omega (\cc )$, the pair $(\cc ,\psi )\in\IBr (G)\times\IBr (\No G(Q))$ satisfies Condition~\ref{Coho}. 
  
If moreover $G$ is perfect or abelian simple,  one says $G$ satisfies ($\ha{\rm iAw}$) to mean that its universal covering $\ha G$ satisfies (iAw).}
\end{defi} 

Navarro-Tiep main reduction theorem for Alperin's weight conjecture (see [NvTi] Theorem~A) can be reformulated as follows :

 \begin{theo}\label{NT}  {\sl If a finite group $G$ is such that any simple subquotient of $G$ satisfies ($\ha{\rm iAw}$), then $G$ satisfies the Alperin weight conjecture. 
 }
\end{theo} 

\begin{rema}\label{RemCoh} {\rm (i) Condition~\ref{Coho} is equivalent to the cohomological condition of [NvTi]~\S 3.3 by the same argument as in [S1]~2.8 which applies over any algebraically closed field $\F$ (not just $\C$). Concerning the finiteness of the cocycles involved (assumed in [S1]~\S 2) see Remark~\ref{CohoRk}.(ii) above. 

(ii) Concerning ($\ha{\rm iAw}$) for a non-abelian simple group $S$, note that it is equivalent to $\ha S/\Ze{\ha S}_p$ satisfying (iAw). This is because all linear representations involved have $\Ze {\ha S}_p ={\rm O}_p(\ha S)$ in their kernel, and one may easily prove a converse of \Pr{CohQu} when $Z'$ is $\Aut (X)$-stable.
} 
\end{rema}

According to [NvTi]~3.2, one of the main steps of the proof of the reduction Theorem is as follows

 \begin{theo}\label{NT32}  {\sl Assume $G$ and $X$ are finite groups with $X\tr G$ and $[G,\Ze X]=\{ 1\}$. Assume moreover $X/\Ze X=S^n$ with $n\geq 1$ and $S$ a simple group satisfying \rm ($\ha{\rm iAw}$). \sl

 Then there exists a bijection $$\Omega\colon \IBr (X)\to \AWe (X) $$
 which is $G$-equivariant, preserves the partition along elements of $\IBr (\Ze X)$ and satisfying, for any $\cc\in\IBr (X)$, $(Q,\psi)\in\Omega (\cc)$, $$|\IBr (G |\cc )|=|\IBr (\No G(Q) |\psi)|.$$}
\end{theo}

\noindent{\it Proof of~\Th{NT32}.} The case where $S$ is abelian (of prime order) is trivial since $X$ is then nilpotent and its only $p$-radical subgroup is normal in $G$. So we assume $S$ simple non-abelian. The $p$-radical subgroups of a direct product are direct products of $p$-radical subgroups of each term, so that $\AWe (\ha S^n)=\AWe (\ha S)^n$ just like we classically have $\IBr (\ha S^n)=\IBr (\ha S)^n$. So, from the fact that $\ha S$ satisfies {\rm (iAw)} with some bijection $\Omega\colon\IBr (\ha S)\to\AWe (\ha S)$, we get a bijection $$\Omega^n\colon\IBr (\ha S^n)\to\AWe (\ha S^n).$$ It is $\Aut (\ha S^n)$-equivariant since $\Om$ is $\Aut (\ha S)$-equivariant and $\Aut (\ha S^n)=\Aut (\ha S){\wr}\Sy_n$ (see the argument in the proof of \Pr{CohPr}). Note that the $\Aut (\ha S)$-equivariance of $\Om$ also ensures the property (*) in \Pr{CohPr}. The cohomological condition for $\Om^n$ (see Condition~\ref{Coho}) is now inherited from the one of $\Om$ thanks to \Pr{CohPr}. 

Further, we have $[X,X]=\ha S^n/Z'$ for some $Z'\inn\Ze {\ha S^n}$. We have $\IBr (\ha S^n/Z')\inn \IBr (\ha S^n)$ and $\AWe (\ha S^n/Z')\inn \AWe (\ha S^n)$ both subsets corresponding with elements of $\IBr (\Ze {\ha S^n})$ with $Z'$ in their kernel. So we deduce from our first step a bijection $$\Om '\colon \IBr ([X,X])\to\AWe ([X,X])$$ which satisfies the cohomological condition thanks to \Pr{CohQu}. The same idea allows to extend it into a bijection $$\Om\colon\IBr (X)\to\AWe (X)$$ since $X$ is a central quotient of $[X,X]\times \Ze X$. Now, $\Om$ preserves the partition along elements of $\IBr (\Ze X)$ and $\Om$ is $G$-equivariant thanks to the $\Aut ([X,X])$-equivariance of $\Om '$ and $[G,\Ze X]=\{ 1\}$. Let $\cc\in\IBr (X)$, $(Q,\psi)\in\Om (\cc)$. The hypotheses of \Th{MainT} are satisfied with $Y:=\No X(Q)$. We have $\No X(Y)=Y$ by $p$-radicality of $Q$, and the other hypotheses in \Th{MainT}.(ii) come from the $G$-equivariance of $\Om$. Then the conclusion of \Th{MainT} gives our last claim.

\qed

\section{Fixed point modules.}

 Here is a case ensuring Condition~\ref{Coho}. We keep $p$ and $\F$.

 \begin{prop}\label{MQ}  {\sl Let $X=[X,X]$ a finite perfect group and $Q\inn X$ a radical $p$-subgroup. Let $V$ be a simple $\F X$-module affording $\cc\in\IBr (X)$, such that 
 
 (i) the fixed point module $V^Q:=\{ v\in V\mid\ q.v=v$ for all $q\in Q\}$ is simple as $\F\No X(Q)$-module, affording $\psi\in\IBr (\No X(Q))$
 
 (ii) $\Aut (X)_{Q,\cc }=\Aut (X)_{Q,\psi}$.
 
  Then $(\cc ,\psi )$ satisfies Condition~\ref{Coho} with respect to $(X,\No X(Q))$.
 }
\end{prop} 

\noindent{\it Proof of~\Pr{MQ}.} Since $X$ is perfect, we have $\Ze {X/\Ze X_p}=\Ze X/\Ze X_p$, while $\Ze X_p$ is in the kernel of all modules involved. So we assume $\Ze X$ is $p'$. Then $Q$ injects in $X/\Ze X$ and we can arrange that the section $\rrr\colon X/\Ze X\to X$ is the identity (hence a group morphism) on $Q$. We denote by $\mathcal{L} $, resp. $\mathcal{M}  $ the irreducible linear representation of $X$, resp. $\mathrm{N}_X(Q)$, on $V$, resp. $V^Q$, so that $\MMM (x)(v)=\LLL (x)(v)$ for any $x\in \No X(Q)$, $v\in V^Q$.

Note that $\Aut (X)_\cc$ transforms $\LLL\circ \rrr$ into equivalent projective representations of $X$.
So, \Pr{ExtPro}.(i) applied to $X/\Ze X\tr\Aut (X)_\cc$ gives us the existence of a projective representation $\mathcal{P}\colon \mathrm{Aut} (X)_\cc\to \mathrm{GL} (V)$ extending the projective representation $\mathcal{L} \circ\rrr$.
Let us show that $$\mathcal{P} (x)(V^Q)\subseteq V^Q\leqno(\mathrm{I})$$ for any $x\in \mathrm{Aut} (X)_{\cc ,Q}$ (note that, $Q$ being $p$-radical, $\Aut (X)_Q=\Aut (X)_{\No X(Q)}$). 

The definition of $V^Q$ and $\rrr$ imply that $$V^Q=\{ v\in V\mid \mathcal{P}(q)(v)=v\text{ for all } q\in Q\}.$$ This also equals $\{ v\in V\mid \mathcal{P}(q)(v)\in\Bbb{F}^\times v$ for all $q\in Q\}$ since $Q$ is a $p$-group on which $\mathcal{P}$ is a group morphism. Now (I) is a consequence of $\mathcal{P} (q)\mathcal{P} (x)\in\Bbb{F}^\times \mathcal{P} (x)\mathcal{P} (x^\mm qx)$ whenever $q\in Q$ and $x\in \mathrm{Aut} (X)_{Q}$. 

Now (I) above allows to define a projective representation $\mathcal{Q} $ of $\mathrm{Aut} (X)_{\cc ,Q}=\mathrm{Aut} (X)_{\No X(Q),\psi }$ (by (ii)) corresponding to the sub-representation of the restriction of $\mathcal{P}$ to $\mathrm{Aut} (X)_{\cc ,Q}$ with space $V^Q$. Then the restriction of $\mathcal{Q}$ to $\No X(Q)/\Ze X$ is $\mathcal{M} \circ\rrr_{|\mathrm{N}_X(Q)/\Ze X}$, and since $V^Q\not= 0$, $\Res^{\Aut (X)_\cc}_{\Aut (X)_{Q,\cc}}\mathcal{P} $ and $\mathcal{Q} $ have same cocycle.

\qed

\begin{rema}\label{MQRk} {\rm  If necessary, it is easy to prove a strengthened version of the above where the fixed point subspace $V^Q$ is replaced by any non-zero $\F\No X(Q)$-module of type $V_i\cap V'_j/V_{i+1}\cap V'_j$ or $V'_i\cap V_j/V'_{i+1}\cap V_j$ where $V_i=J^iV$ and $V'_i=\{v\in V\mid J^{n-i}.v=0\}$ for $J$ the Jacobson radical of $\F Q$, $n$ the nilpotence index of $J$, and $0\leq i\leq n-1$.  }\end{rema}

We show how the above proposition allows a quite elementary checking of {\rm (iAw)} for quasi-simple groups of Lie type of characteristic $p$, thus giving another proof of [NvTi] Theorem~C.

 \begin{prop}\label{BNdef}  {\sl Let $G$ be a finite perfect group endowed with a strongly split BN-pair of characteristic $p$ (see [CaEn]~2.20). Then $G$ satisfies {\rm (iAw)} for the prime $p$.
 }
\end{prop} 

\noindent{\it Proof of~Proposition~\ref{BNdef}.} Let $Y:=\F [G/U]$ where $U$ is a Sylow $p$-subgroup of $G$. Recall the so-called Green correspondence $$M\mapsto f(M)$$ which associates to any indecomposable $\F G$-module $M$ with vertex $Q$, an indecomposable $\No G(Q)$-module $f(M)$ with same vertex, see [B]~\S 3.12, [NgTs]~\S 4.4. By a general fact (see [Al]~Lemma~1), when $(Q,S)$ is a pair with $Q$ a $p$-subgroup of $G$ and $S$ a simple $\F\No G(Q)$-module affording an element of $\IBr_Q(\No G(Q))$, the Green vertex of $S$ is $Q$ and the Green correspondent $f^\mm(S)$ of $S$ is a direct summand of $Y$. Denote $$f'(S):=f^\mm(S)/J(\F G).f^\mm(S)$$ the largest semi-simple quotient of $f^\mm(S)$ as $\F G$-module. Recall (AWC) from the preceding section. Let us show that in our case 

(a) the map $$f'(S)\mapsfrom (Q,S)$$ is a bijection realizing (AWC)  \and 

(b) $S\cong f'(S)^Q$ (fixed points). 

This will complete our proof since the above map is clearly $\Aut (G)$-equivariant, preserves characters of $\Ze G_\pp$ (or even $p$-blocks), and \Pr{MQ} gives the required cohomological condition (note that hypothesis (ii) in that proposition is ensured by equivariance). 

To prove (a) and (b), we follow the checking of (AWC) for those groups given in [Ca] (see also [CaEn]~\S 6.3). Let us recall the labelling of $\IBr (G)$ and of indecomposable summands of $Y$ due to Green-Sawada-Tinberg (see [CuRe]~\S 72.B, [CaEn]~\S 6, [T]). The axioms of BN-pairs ensure that $G$ has subgroups $B$ (Borel), $T$ (torus) and an associated set $\Pi$ (fundamental roots) such that the subgroups of $G$ containing $B$ (parabolic subgroups) are parametrized $\DD\mapsto P_\DD$ by subsets $\Delta$ of $\Pi$. The BN-pair being strongly split, each $P_\Delta$ is a semi-direct product $P_\DD =\No G(U_\DD)=U_\DD\smd L_\DD$ where $U_\DD={\rm O}_p(P_\DD)$ and $L_\DD$ (Levi subgroup) is a group of the same type as $G$ with same torus $T$, Borel subgroup $B\cap L_\DD$ and $\DD$ as set of fundamental roots. One also has subgroups $T_\DD\inn T$ in (increasing) bijection with subsets of $\Pi$. 

One has $$Y=\oplus_{(\la ,\DD)}Y(\la ,\DD )$$ where the sum is over ``admissible pairs", i.e. pairs $(\la ,\DD)$ where $\la\in\IBr (T)$ and $\DD\inn\Pi$ is such that $\la (T_\DD)=\{ 1\}$ (see [T]~\S 2, [Ca]~\S B.9.2). Each $Y(\la ,\DD)$ is indecomposable with head $M_G(\la ,\DD)$ a simple $\F G$-module and $$(\la ,\DD)\mapsto M_G(\la ,\DD )$$ provides a parametrization of $\IBr (G)$ by admissible pairs (see [CaEn]~6.12 or [CuRe]~72.28). Moreover $Y(\la ,\DD)$ has Green vertex $U_\DD$ and Green correspondent $M_{L_\DD}(\la ,\DD)$ as $\F(\No G(U_\DD )/U_\DD)$-module (see [Ca]~Proposition~6, [T]~\S 3). This gives (a), both sets involved being in bijection with admissible pairs. 

Now (b) rewrites $$M_{L_\DD}(\la ,\DD )=M_G(\la ,\DD)^{U_\DD},$$ which is a special case of [CaEn]~6.12.(iii) (see also [GLS]~2.8.11). 

\qed

Except in three cases, finite simple groups $S$ of Lie type of characteristic $p$ have a universal covering $\ha S$ such that $\ha S/\Ze{\ha S}_p$ is also a group of Lie type endowed with a strongly split BN-pair of same characteristic (see the tables in [GLS]~\S 6.1]). Then the above gives another proof of [NvTi]~Theorem~C. 

The three exceptions are for the prime $p=2$.  The simple groups $\PSL_2(9)$, $\PSU_3(3)$ are of characteristic $3$ but sometimes also considered as groups of characteristic $2$ since isomorphic to $[\Sp_4(2),\Sp_4(2)] $, $[G_2(2),G_2(2)] $ respectively. For those two and for $[^2F_4(2),^2F_4(2)]$, a direct checking using available tables of 2-modular characters has been made, see [NvTi]~6.1.

\section*{References.}

[Al] J.L. Alperin, Weights for finite groups, \sl in Proc. Symp. pure Math., \bf 47 \rm I (1987), 369-379 .

[As] M. Aschbacher, \it Finite Group Theory, \rm Cambridge University Press, Cambridge, 1986.

[B] D. Benson, \it Representations and Cohomology 
I: Basic Representation Theory of Finite 
Groups and Associative Algebras, \rm Cambridge University Press, Cambridge, 1991.


[Ca] M. Cabanes, Brauer morphisms between modular Hecke algebras, \sl J. Algebra, \bf 115-1\rm (1988), 1-31.

[CaEn] M. Cabanes and M. Enguehard, {\it Representation theory of finite reductive groups}, Cambridge University Press, Cambridge, 2004.

[CuRe] C.W. Curtis and I. Reiner, {\it Methods of representation theory with applications to finite groups and orders}, Wiley, New York, 1981.


[GLS] D. Gorenstein, R. Lyons and R. Solomon, {\it The classification of the finite
simple groups,
Number 3.} {\sl Mathematical Surveys and Monographs, Amer. Math. Soc.}, Providence, 1998.


[IMN] M. Isaacs, G. Malle, and G. Navarro, A reduction theorem for the
McKay conjecture, \sl Invent. Math., \bf 170 \rm (2007),
33--101.

[NgTs] H. Nagao and Y. Tsushima, \it Representations of Finite Groups, \rm Academic, Boston, 1989.

[Nv] G. Navarro, Characters and Blocks of Finite Groups, Cambridge University Press, Cambridge, 1998.

[NvTi] G. Navarro, P.H. Tiep, A reduction theorem for the Alperin weight conjecture, \sl Invent. Math., \bf 184 \rm (2011),
529--565.

[P] L. Puig, On the reduction of Alperin's conjecture to the quasi-simple groups, 
\sl J. Algebra, \bf 328 \rm (2011), 372--398.

{[S1]}  B. Sp\"ath, Inductive McKay condition in defining characteristic, \sl Bull. London Math. Soc., \bf 44 \rm (2012), 426--438. \rm

{[S2]}  B. Sp\"ath, A reduction theorem for the Alperin-McKay conjecture, \sl J. reine angew. Math. \rm (2012), \it to appear.\rm

{[S3]}  B. Sp\"ath,
 A reduction theorem for the blockwise {A}lperin weight conjecture, \sl J. Group Theory, \rm (2012),
\it to appear, \tt dx.doi.org/10.1515/jgt-2012-0032.\rm

[T] N.B. Tinberg,  Some indecomposable modules of groups with split $(B,\,N)$-pairs, \sl J. Algebra, \bf 61 \rm (1979), no. 2, 508--526

\bigskip

{}

\bigskip

{}

\bigskip

{}

\rightline{\small Marc Cabanes, CNRS-Universit\'e Paris 7, Paris, France}

\rightline{\small {cabanes@math.jussieu.fr}}

\rightline{\small {\tt http://people.math.jussieu.fr/\~{}cabanes/}}

\end{document}